\documentclass[12pt]{article}

\usepackage{amssymb,amsmath,amsthm}
\usepackage{amsfonts}
\usepackage{enumerate}

\usepackage{fancyhdr}	
\usepackage{mdframed}
\usepackage{geometry}
\newcommand\blfootnote[1]{%
  \begingroup
  \renewcommand\thefootnote{}\footnote{#1}%
  \addtocounter{footnote}{-1}%
  \endgroup
}

\newtheoremstyle{mytheoremstyle} 
    {10pt}                    
    {10pt}                    
    {\normalfont}                   
    {}                           
    {\bfseries}                   
    {.}                          
    {0.3cm}                       
    {}  
\theoremstyle{mytheoremstyle}

\theoremstyle{plain}

\fancypagestyle{plain}{%
  \fancyhf{}%
}

\newtheorem{theorem}{Theorem}[section]

\newtheorem{lemma}[theorem]{Lemma}
\newtheorem{proposition}[theorem]{Proposition}

\newtheorem{remark}[theorem]{Remark}

\newtheorem{corollary}[theorem]{Corollary}
\newtheorem{problem}[theorem]{Problem}

\begin{document}

\title{ $k$-spectrally monomorphic tournaments}
\author{
	Abderrahim Boussa\"{\i}ri, Imane Souktani, Imane Talbaoui, Mohamed Zouagui
}

\maketitle

\begin{abstract}
 A tournament is $k$-spectrally monomorphic if all the $k\times k$
principal submatrices of its adjacency matrix have the same characteristic
polynomial. Transitive $n$-tournaments are trivially $k$-spectrally monomorphic. We show that there are no other for
$k\in  \{  3,\ldots,n-3 \}  $. Furthermore, we prove that for $n\geq
5$, a non-transitive $n$-tournament is $(n-2)$-spectrally monomorphic if and
only if it is doubly regular. Finally, we give some results on
 $(n-1)$-spectrally monomorphic regular tournaments.
\end{abstract}



\section{Introduction}
\blfootnote{Corresponding author: Abderrahim Boussa\"{\i}ri. Email: aboussairi@hotmail.com}
\blfootnote{Laboratoire Topologie, Alg\`ebre, G\'eom\'etrie et Math\'ematiques Discr\`etes, Facult\'e des Sciences A\"in Chock, Hassan II University of Casablanca, Maroc.}						
An $n$-\emph{tournament} $T$ is a digraph with $n$ vertices in which every
pair of vertices is jointed by exactly one arc. If the arc joining two given vertices
$u$ and $v$ of $T$ is directed from $u$ to $v$, we say that $u$ \emph{dominates} $v$.
A tournament $T$ is \emph{transitive} if whenever $u$ dominates $v$ and $v$
dominates $w$, then $u$ dominates $w$. With respect to an ordering
$v_{1},\ldots,v_{n}$ of the vertices of $T$, the \emph{adjacency matrix} of
$T$ is the $n\times n$ matrix $A=(a_{ij})$ in which $a_{ij}$ is $1$ if $v_{i}$
dominates $v_{j}$ and $0$ otherwise. The \emph{skew-adjacency matrix} of $T$ is
$S:=A-A^{\top}$, where $A^{\top}$ is the transpose of $A$. Clearly, 
$A+A^{\top}=J_{n}-I_{n}$, where $I_{n}$ and $J_{n}$ denote respectively the
$n\times n$ identity matrix and the all-ones matrix. The \emph{characteristic
polynomial} and the \emph{skew-characteristic polynomial} of $T$ are
 the characteristic polynomial of $A$ and $S$ respectively.

Following \cite{attas2021characterization}, a square matrix is $k$-\emph{spectrally
monomorphic} if all its $k\times k$ principal submatrices have the same
characteristic polynomial. The main result of \cite{attas2021characterization} is the characterization of
$k$-spectrally monomorphic Hermitian matrices for $k\in \{  3,\ldots
,n-3 \}$.

A tournament is $k$-\emph{spectrally monomorphic} if its adjacency matrix is
$k$-spectrally monomorphic. Spectral monomorphy is a weakening of monomorphy for tournaments. An
$n$-tournament is said to be $k$-\emph{monomorphic} if all its $k$-subtournaments are
isomorphic. For example, transitive tournaments are $k$-monomorphic for
every $k$. Conversely, it follows from a combinatorial lemma of Pouzet
\cite{Pouzet1976application} that for any integer $n\geq6$ and any
$k\in  \{  3,\ldots,n-3 \}  $, $k$-monomorphic $n$-tournaments are
transitive. Using similar techniques, we prove that the $k$-monomorphy can be
relaxed to $k$-spectral monomorphy.

Jean \cite{jean1969line} and
later Pouzet \cite{Pouzet1978} proved that a non-transitive $n$-tournament is
$(n-2)$-monomorphic if and only if it is arc-symmetric, that is, its
automorphism group acts transitively on the set of its arcs. Arc-symmetric
tournaments were characterized independently by Kantor \cite{Kantor1969} and
Berggren \cite{Berggren1972}. They showed that an arc-symmetric tournament is
isomorphic to one of the quadratic residue tournaments, defined on a finite
field $\mathbb{F}_{p^{r}}$, $p^{r}\equiv3\pmod4$ by the
following rule: $x$ dominates $y$ if and only if $x-y$ is a non-zero square in
$\mathbb{F}_{p^{r}}$. In a quadratic residue tournament, all pairs of vertices
jointly dominate the same number of vertices. Tournaments with this property
are called \emph{doubly regular}. There are many structural and spectral
characterizations of doubly regular tournaments \cite{Reid1972, Muller1974,
Zagaglia1982, Rowlinson1986}. In this paper, we obtain the following.

\begin{theorem}
\label{main2}An $n$-tournament with $n\geq5$ is $(n-2)$-spectrally monomorphic
if and only if it is transitive or doubly regular. Moreover, the common
characteristic polynomial of all the $(n-2)$-subtournaments of a doubly regular
$n$-tournament is
\[
\left(z^{2}+z+\tfrac{n+1}{4}\right)^{\tfrac{n-5}{2}}\left(z^{3}-\tfrac{n-5}{2}z^{2}-\tfrac{n-5}%
{2}z-\tfrac{n-3}{4}\right)\mbox{.}%
\]
\end{theorem}

The characterization of $(n-1)$-spectrally monomorphic $n$-tournaments
seems difficult. It is perhaps of comparable difficulty as the open problem of
characterizing $(n-1)$-monomorphic $n$-tournaments \cite[Problem 43 p.
252]{Kotzig-1monomorph}. Doubly regular $n$-tournaments are $(n-1)$-spectrally
monomorphic. Indeed, given a doubly regular $n$-tournament, Theorem 1.1 and
Theorem 2.7 of \cite{Nozaki2012} imply that the characteristic polynomial of
all its $(n-1)$-subtournaments is $\left(z^{2}+z+\frac{n+1}{4}\right)^{\frac{n-3}{2}%
}\left(z^{2}+\frac{3-n}{2}z+\frac{3-n}{4}\right)$. As for regular $n$-tournaments, they
are not always $(n-1)$-spectrally monomorphic.

\begin{theorem}
\label{counter-example}For every integer of the form $n=3^{k}\cdot7$, there
exists a non $(n-1)$-spectrally monomorphic regular $n$-tournament.
\end{theorem}

\section{Properties of $k$-spectrally monomorphic matrices}

Let $M$ be an $n\times n$ matrix. The characteristic polynomial of $M$ is
$P_{M}(z):=\det(zI_{n}-M)$. For a subset $\alpha$\ of $[  n]
:=\{1,\ldots,n\}$, we denote by $M[\alpha]$ the principal submatrix of $M$
whose rows and columns are indexed by the elements of $\alpha$. For $i\in [
n]  $, we denote $M[[  n]  \setminus  \{  i \}  ]$ simply by
$M_{i}$.
Let $P_{M}(z)=z^{n}+a_{1}z^{n-1}+a_{2}z^{n-2}+a_{3}z^{n-3}+\cdots+a_{n}$. It
is well-known (see for example \cite[p. 494]{Meyer2000}), that%

\begin{equation}
a_{k}=(-1)^{k}\sum_{\alpha \in \binom{[n]}{k}}\det
M[\alpha] \label{eq1}%
\end{equation}
for every $k$ in $[  n]  $.

\begin{remark}
If $M$ is the adjacency matrix of a tournament $T$,
then $a_{1}=a_{2}=0$ and $-a_{3}$ is the number of $3$-cycles of $T$.
\end{remark}

The derivative of the characteristic polynomial of $M$ is equal to the sum of
the characteristic polynomials of its $(n-1)\times(n-1)$ principal submatrices
\cite[Theorem 1]{Schwenk79}. Applying this fact iteratively, we get
\begin{equation}
P_{M}^{(n-k)}(z)=(n-k)!\sum_{\alpha \in \binom{[  n]  }{k}}{P_{M[\alpha]}(z)}\mbox{,}\label{schwenk2}%
\end{equation}
where $P_{M}^{(i)}(z)$ denotes the $i$-th derivative of $P_{M}(z)$. As a
consequence, we have

\begin{corollary}
\label{corol-poly-commun}Suppose that $M$ is $k$-spectrally monomorphic and
let $Q$ be the common characteristic polynomial of its $k\times k$ principal
submatrices. Then
\[
P_{M}^{(n-k)}(z)=\frac{n!}{k!}Q(z)\mbox{.}%
\]

\end{corollary}

In \cite{attas2021characterization}, Attas et al. showed that spectrally monomorphic
matrices enjoy a hereditary property.

\begin{proposition}
\label{prop inf kspectra}  $k$-spectrally monomorphic $n\times n$ complex matrices 
are $l$-spectrally monomorphic for each $l\in  \{  1,\ldots
,\min(k,n-k) \}  $.
\end{proposition}

A slight modification of the proof of this proposition yields the following
result, which plays a key role in our study of $k$-spectrally monomorphic tournaments.

\begin{proposition}
\label{conseq 2 pouzet}Let $M$ be an $n\times n$ complex matrix, $k\in  \{
3,\ldots,n-1 \}  $ and $p\leq k$. If $M$ is $k$-spectrally
monomorphic, then for every subset $\beta \subseteq [n]  $ of
size at most $n-k$, the number $\sum_{\alpha \supseteq \beta, \vert
\alpha  \vert =p}{\det(M[\alpha])}$ depends only on the cardinality of
$\beta$.
\end{proposition}

\section{$k$-spectrally monomorphic tournaments for \\ $k\in  \{
3,\ldots,n-2 \}  $}

Obviously, all tournaments are $2$-spectrally monomorphic, so we will
only consider $k$-spectrally monomorphic tournaments for $k\geq3$.

\begin{proposition}
\label{main1}
$k$-spectrally monomorphic $n$-tournaments are transitive for every $n\geq6$ and $k\in  \{  3,\ldots,n-3 \}  $.
\end{proposition}

\begin{proof}
Let $k\in  \{  3,\ldots,n-3 \}  $ and let $T$ be a $k$-spectrally
monomorphic $n$-tournament. As $n\geq 6 $, by Proposition \ref{prop inf kspectra}, $T$ is
$3$-spectrally monomorphic. The two tournaments with three vertices, namely
the $3$-cycle and the transitive tournament, have different characteristic
polynomials. Moreover, every tournament with more than three vertices contains
a transitive $3$-tournament. Hence $T$ is necessarily transitive because, by
hypothesis, all its $3$-subtournaments have the same characteristic polynomial.
\end{proof}

Let $T$ be an $n$-tournament with vertex set $V$ and let $x\in V$. The
\emph{out-neighborhood} of the vertex $x$ is $N_{T}^{+}(x):=\{y\in V: x \mbox{ dominates } y\}$
and its \emph{in-neighborhood} is $N_{T}^{-}(x):=\{y\in V: y \mbox{ dominates } x\}$.
The \emph{out-degree} (resp. the \emph{in-degree}) of $x$ is $\delta_{T}%
^{+}(x):= \vert N_{T}^{+}(x) \vert $ (resp. $\delta_{T}%
^{-}(x):= \vert N_{T}^{-}(x) \vert $). The tournament $T$ is
\emph{regular} if there exists $k\geq1$ such that the out-degree of all its
vertices is $k$. In this case, $n$ is odd and $k=\tfrac{n-1}{2}$. The
tournament $T$ is \emph{near-regular} if $n$ is even and the out-degree of
every vertex is $\tfrac{n-2}{2}$ or $\tfrac{n}{2}$.

Let $x$ and $y$ be distinct vertices of $T$. The \emph{out-degree} (resp. the
\emph{in-degree}) of $(  x,y)  $ is $\delta_{T}^{+}%
(x,y):= \vert N_{T}^{+}(x)\cap N_{T}^{+}(y) \vert $ (resp.
$\delta_{T}^{-}(x,y):= \vert N_{T}^{-}(x)\cap N_{T}^{-}(y) \vert $).
With these notations, the tournament $T$ is doubly regular if the out-degree
of all pairs of vertices is the same.

Recall the following fundamental properties of doubly regular tournaments (for
the proof, see \cite{Muller1974}).

\begin{proposition}
\label{DoublyRegular} Let $T$ be a doubly regular $n$-tournament. Then $T$ is
regular and there exists $t\geq0$ such that $n=4t+3$. Moreover, if $x$ and $y$
are two distinct vertices of $T$ such that $x\ $dominates $y$ then
\[
 \vert N_{T}^{+}(x)\cap N_{T}^{+}(y) \vert = \vert N_{T}%
^{-}(x)\cap N_{T}^{-}(y) \vert = \vert N_{T}^{+}(x)\cap N_{T}%
^{-}(y) \vert =t
\]%
\[
\mbox{and } \vert N_{T}^{-}(x)\cap N_{T}^{+}(y) \vert =t+1.
\]

\end{proposition}

Notice that $ \vert N_{T}^{-}(x)\cap N_{T}^{+}(y) \vert $ is the
number of $3$-cycles containing $x$ and $y$. A tournament in which the number
of $3$-cycles containing any pair of vertices is a positive constant, say $k$,
is called \emph{homogeneous}. Kotzig \cite{Kotzig1969} proved that such a
tournament has $4k-1$ vertices. Reid and Brown \cite{Reid1972} proved that
homogeneous tournaments are doubly regular. Consequently

\begin{proposition}
\label{transitive or doubly}Let $T$ be a tournament. The following statements
are equivalent
\begin{enumerate}
\item[i)] Every pair of vertices is contained in the same number of $3$-cycles.
\item[ii)] $T$ is transitive or doubly regular.
\end{enumerate}
\end{proposition}

Throughout this paper, the zero and the all-ones $p\times q$ matrices are denoted by
$O_{p,q}$ and $J_{p,q}$. The matrices $O_{p,p}$ and $J_{p,p}$ are simply
written as $O_{p}$ and $J_{p}$. The zero and the all-ones column vectors of
size $p$ are denoted respectively by $\mathbf{0}_{p}$ and $\mathbf{1}_{p}$.

Let $A$ be the adjacency matrix of an $n$-tournament $T$ with respect to an
ordering $x_{1},\ldots,x_{n}$ of its vertices. For any $i,j\in \{1,\ldots
,n\}$,\ the $(i,j)$-entry of $AA^{\top}$ is

%

\begin{equation}\label{AAt}%
\begin{cases}
\delta_{T}^{+}(x_{i})&\mbox{ if }i=j\mbox{,}\\
\delta_{T}^{+}(x_{i},x_{j})&\mbox{ if }i\neq j \mbox{.}
\end{cases}  
\end{equation}

If $T$ is doubly regular with $n=4t+3$ vertices then $AA^{\top}=tJ_{n}%
+(t+1)I_{n}$. This equality was used in \cite[Proposition 3.1]{Caen1992} to
prove that%

\begin{equation}
P_{A}(z)=\left(  z-\tfrac{n-1}{2}\right)  \left(  z^{2}+z+\tfrac{n+1}{4}\right)
^{\tfrac{n-1}{2}}\mbox{.} \label{chact doubly}%
\end{equation}

\begin{proof}[Proof of Theorem \ref{main2}]
Let $T$ be an $(n-2)$-spectrally monomorphic $n$-tournament with vertex set $ \{  v_{1},\ldots,v_{n} \}  $ and
adjacency matrix $A$. Let $v_{i}$ and $v_{j}$ be two distinct vertices of $T$. From
Proposition \ref{conseq 2 pouzet}, the sum $\sum_{l\in [  n]
\setminus  \{  i,j \} }{\det(A[ \{  i,j,l \}  ])}$,
which is exactly the number of $3$-cycles containing $v_{i}$ and $v_{j}$, does
not depend on $v_{i}$ and $v_{j}$. Then, by Proposition
\ref{transitive or doubly}, $T$ is doubly regular or transitive.

Conversely, let $T$ be a doubly regular tournament with $n=4t+3$ vertices and
let $x$ and $y$ be two distinct vertices of $T$ such that $x$ dominates $y$.
Denote by $R$ the tournament obtained from $T$ by removing $x$ and $y$. The
vertex set of $R$ is partitioned into the following four subsets: $N_{T}%
^{+}(x)\cap N_{T}^{+}(y)$, $N_{T}^{-}(x)\cap N_{T}^{+}(y)$, $N_{T}^{+}(x)\cap
N_{T}^{-}(y)$ and $N_{T}^{-}(x)\cap N_{T}^{-}(y)$. Let $B$ be the adjacency
matrix of $R$. Using \eqref{AAt} and Proposition \ref{DoublyRegular}, we get
\[
BB^{\top}-(t+1)I_{4t+1}=
\begin{pmatrix}
tJ_{t} & tJ_{t,t+1} & tJ_{t} & tJ_{t}\\
tJ_{t+1,t} & (t-1)J_{t+1} & tJ_{t+1,t} & (t-1)J_{t+1,t}\\
tJ_{t} & tJ_{t,t+1} & (t-1)J_{t} & (t-1)J_{t}\\
tJ_{t} & (t-1)J_{t,t+1} & (t-1)J_{t} & (t-2)J_{t}%
\end{pmatrix}
  \mbox{.}%
\]
As in the proof of Proposition 3.1 of \cite{Caen1992}, we have
\begin{align*}
P_{B}^{2}(z)  &  =\det\left(zI_{4t+1}-B\right)\det\left(zI_{4t+1}-B^{\top}\right)\\
&  =\det\left(z^{2}I_{4t+1}-z\left(B+B^{\top}\right)+BB^{\top}\right)\\
&  =\det\left(z^{2}I_{4t+1}-z\left(J_{4t+1}-I_{4t+1}\right)+BB^{\top}\right)\mbox{.}%
\end{align*}
Then $P_{B}^{2}(z)$, and a fortiori $P_{B}(z)$, do not depend on the choice of
$x$ and $y$. Hence $T$ is $(n-2)$-spectrally monomorphic.

Let $Q$ be the common characteristic polynomial of the $(n-2)$-subtournaments
of $T$. Using \eqref{chact doubly} and Corollary \ref{corol-poly-commun}, we
get
\begin{align*}
Q(z)  &  =\tfrac{1}{n(n-1)}P_{A}^{(2)}(z)\\
&  =\left(  z^{2}+z+\tfrac{n+1}{4}\right)  ^{\tfrac{n-5}{2}}\left(  z^{3}%
-\tfrac{n-5}{2}z^{2}-\tfrac{n-5}{2}z-\tfrac{n-3}{4}\right)\mbox{.}%
\end{align*}

\end{proof}

\section{A family of regular $n$-tournaments which are not $(n-1)$-spectrally monomorphic}

As mentioned in the introduction, regular $n$-tournaments are not always
$(n-1)$-spectrally monomorphic. The smallest counter-example has $7$
vertices, its adjacency matrix is%
\[
A=
\begin{pmatrix}
0 & 0 & 0 & 0 & 1 & 1 & 1\\
1 & 0 & 0 & 1 & 0 & 0 & 1\\
1 & 1 & 0 & 0 & 1 & 0 & 0\\
1 & 0 & 1 & 0 & 0 & 1 & 0\\
0 & 1 & 0 & 1 & 0 & 0 & 1\\
0 & 1 & 1 & 0 & 1 & 0 & 0\\
0 & 0 & 1 & 1 & 0 & 1 & 0
\end{pmatrix} \mbox{.}%
\]

To obtain an infinite family of counter-examples, we use the following
construction. Let $T_{1}$, $T_{2}$\ and $T_{3}$\ be three regular
$n$-tournaments with disjoint vertex sets $V_{1}= \{  v_{1},\ldots
,v_{n} \}  $, $V_{2}= \{  v_{n+1},\ldots,v_{2n} \}  $ and
$V_{3}= \{  v_{2n+1},\ldots,v_{3n} \}  $ respectively. Consider the
$3n$-tournament $T$ with vertex set $V=V_{1}\cup V_{2}\cup V_{3}$, obtained
from $T_{1}$, $T_{2}$\ and $T_{3}$ by adding arcs from $V_{1}$ to $V_{2}$,
$V_{2}$ to $V_{3}$ and $V_{3}$ to $V_{1}$.

\begin{proposition}
\label{construction counter}The $3n$-tournament $T$ is regular. Moreover, if
at least one of the three tournaments $T_{1}$, $T_{2}$\ and $T_{3}$\ is not
$(n-1)$-spectrally monomorphic then $T$ is not $(3n-1)$%
-spectrally monomorphic.
\end{proposition}

The proof of this proposition requires the following lemma.

\begin{lemma}
\label{regular exp tournament}Let $T$ be a regular $n$-tournament with
adjacency matrix $A$. Then, for every $i\in [  n]  $, we have%
\[
P_{A}(z)=\tfrac{1}{2}\left( z-\tfrac{n-1}{2}\right)  [  \left(n+2z+1\right) P_{A_{i}%
}(z)+\left(n-2z-1\right)P_{A_{i}}(-z-1)]  \mbox{.}%
\]

\end{lemma}

\begin{proof}
 Without loss of
generality, we assume that $i=1$. The matrix $A$ is partitioned as follows%
\[
A=
\begin{pmatrix}
0 & w^{\top}\\
\mathbf{1}_{n-1}-w & A_{1}%
\end{pmatrix}\mbox{,}
\]
where $w$ is an $(n-1)$-column vector with coordinates in $ \{
0,1 \} $.
Since $T$ is regular, $n=2k+1$ for some positive integer $k$, moreover, $A\mathbf{1}_{n}=k\mathbf{1}_{n}$ and $\mathbf{1}%
_{n}^{\top}A=k\mathbf{1}_{n}^{\top}$.
Then
\begin{align}
w^{\top}\mathbf{1}_{n-1} & =\mathbf{1}_{n-1}^{\top}w=k\mbox{;}\\
A_{1}\mathbf{1}_{n-1} & =(k-1)\mathbf{1}_{n-1}+w\mbox{;}\\
\mathbf{1}_{n-1}^{\top}A_{1} & =k\mathbf{1}_{n-1}^{\top}-w^{t}\mbox{.}%
\end{align}
We have%
\begin{align*}
P_{A}(z) &  =\det\left(zI_{n}-A\right)\\
&  = 
\begin{vmatrix}
z & -w^{\top}\\
w-\mathbf{1}_{n-1} & zI_{n-1}-A_{1}%
\end{vmatrix}\mbox{.}
\end{align*}
By adding the sum of the last $(n-1)$ columns of the matrix $zI_{n}-A$ to its
first column, we get%
\begin{align*}
P_{A}(z) &  =
\begin{vmatrix}
z-w^{\top}\mathbf{1}_{n-1} & -w^{\top}\\
w-\mathbf{1}_{n-1}+\left(zI_{n-1}-A_{1}\right)\mathbf{1}_{n-1} & zI_{n-1}-A_{1}%
\end{vmatrix}
  \\
&  = 
\begin{vmatrix}
z-k & -w^{\top}\\
\left(z-k\right)\mathbf{1}_{n-1} & zI_{n-1}-A_{1}%
\end{vmatrix}
  \mbox{.}%
\end{align*}
Adding the sum of the last $(n-1)$ rows of the above matrix to its first row,
we get
\begin{align*}
P_{A}(z) &  =
\begin{vmatrix}
z-k+\mathbf{1}_{n-1}^{\top}\left(z-k\right)\mathbf{1}_{n-1} & -w^{\top}+\mathbf{1}%
_{n-1}^{\top}\left(  zI_{n-1}-A_{1}\right)  \\
\left(z-k\right)\mathbf{1}_{n-1} & zI_{n-1}-A_{1}%
\end{vmatrix}
  \\
&  = 
\begin{vmatrix}
n\left(z-k\right) & \left(  z-k\right)  \mathbf{1}_{n-1}^{\top}\\
\left(z-k\right)\mathbf{1}_{n-1} & zI_{n-1}-A_{1}%
\end{vmatrix}
 \\
&  =n 
\begin{vmatrix}
\left(z-k\right) & \tfrac{1}{n}\left(  z-k\right)  \mathbf{1}_{n-1}^{\top}\\
\left(z-k\right)\mathbf{1}_{n-1} & zI_{n-1}-A_{1}%
\end{vmatrix}
  \mbox{.}%
\end{align*}
Now, subtracting the first row from each of the last $n-1$ rows, we obtain%
\begin{align*}
P_{A}(z) &  =n
\begin{vmatrix}
\left(z-k\right) & \tfrac{1}{n}\left(  z-k\right)  \mathbf{1}_{n-1}^{\top}\\
\left(z-k\right)\mathbf{1}_{n-1} & zI_{n-1}-A_{1}%
\end{vmatrix}
  \\
&  =n
\begin{vmatrix}
\left(z-k\right) & \tfrac{1}{n}\left(  z-k\right)  \mathbf{1}_{n-1}^{\top}\\
\mathbf{0}_{n-1} & zI_{n-1}-A_{1}-\tfrac{1}{n}\left(  z-k\right)  J_{n-1}%
\end{vmatrix}
  \\
&  =n\left(z-k\right)\det\left(zI_{n-1}-A_{1}-\tfrac{1}{n}\left(z-k\right)J_{n-1}\right)\mbox{.}%
\end{align*}
We conclude by the first assertion of lemma below.
\end{proof}

\begin{lemma}
\label{poly-with indetr}\cite{Gregory1993} Let $T$ be an $n$-tournament with
adjacency matrix $A$. For any scalar $\lambda$, we have

\begin{enumerate}
\item[i)] The characteristic polynomial of $A+\lambda J_{n}$ is
\[
P_{A+\lambda J_{n}}(z)=(\lambda+1)P_{A}(z)-(-1)^{n}\lambda P_{A}(-z-1)\mbox{.}%
\]
\item[ii)] If $T$ is regular, then
\[
\left(  z-\tfrac{n-1}{2}\right)  P_{A+\lambda J_{n}}(z)=\left(  z-n\lambda
-\tfrac{n-1}{2}\right)  P_{A}(z)\mbox{.}%
\]

\end{enumerate}
\end{lemma}

\begin{remark}
Lemma \ref{regular exp tournament} shows that the characteristic
polynomial of a regular $n$-tournament is determined by the characteristic
polynomial of any of its $(n-1)$-subtournaments. The technique we used in the
proof is borrowed from that of \cite[Theorem 1.1]{Nozaki2012}.
\end{remark}

Now, we are able to prove Proposition \ref{construction counter}.

\begin{proof}[Proof of Proposition \ref{construction counter}]
 Obviously, the tournament $T$ is
regular. Let $A$, $B$ and $C$ be the adjacency matrices of $T_{1}$,
$T_{2}$\ and $T_{3}$ respectively. The adjacency matrix of $T$ is
\[
M=
\begin{pmatrix}
A & J_{n} & O_{n}\\
O_{n} & B & J_{n}\\
J_{n} & O_{n} & C%
\end{pmatrix}\mbox{.}%
\]

Without loss of generality, suppose that $T_{3}$ is not $(n-1)$%
-spectrally monomorphic and let $i\in [n]  $. The adjacency
matrix of $T-v_{i+2n}$ is
\[
M_{i+2n}=
\begin{pmatrix}
A & J_{n} & O_{n,n-1}\\
O_{n} & B & J_{n,n-1}\\
J_{n-1,n} & O_{n-1,n} & C_{i}%
\end{pmatrix}\mbox{,}%
\]
where $C_{i}$ is the adjacency matrix of $T_{3}-v_{i+2n}$.
The characteristic polynomial of $M_{i+2n}$ is%
\[
P_{M_{i+2n}}(z)= 
\begin{vmatrix}
zI_{n}-A & -J_{n} & O_{n,n-1}\\
O_{n} & zI_{n}-B & -J_{n,n-1}\\
-J_{n-1,n} & O_{n-1,n} & zI_{n-1}-C_{i}%
\end{vmatrix}
\mbox{.}%
\]
Multiplying the last $n-1$ rows by $(z-\tfrac{n-1}{2})$, we get%

\[
P_{M_{i+2n}}(z)
=\tfrac{1}{\left(z-\tfrac{n-1}{2}\right)^{n-1}}
\begin{vmatrix}
zI_{n}-A & -J_{n} & O_{n,n-1}\\
O_{n} & zI_{n}-B & -J_{n,n-1}\\
-\left(z-\tfrac{n-1}{2}\right)J_{n-1,n} & O_{n-1,n} & \left(z-\tfrac{n-1}{2}\right)\left(zI_{n-1}
-C_{i}\right)
\end{vmatrix}\mbox{.}
 \] 
As $T_{1}$ is regular, the sum of the first $n$ rows is
\[
\left(z-\tfrac{n-1}{2},\ldots,z-\tfrac{n-1}{2},-n,\ldots,-n,0,\ldots,0\right)\mbox{.}%
\]
Adding this sum to each of the last $n-1$ rows, we obtain%
\begin{align*}
P_{M_{i+2n}}(z)  &  =\tfrac{1}{\left(z-\tfrac{n-1}{2}\right)^{n-1}} 
\begin{vmatrix}
zI_{n}-A & -J_{n} & O_{n,n-1}\\
O_{n} & zI_{n}-B & -J_{n,n-1}\\
O_{n-1,n} & -nJ_{n-1,n} & \left(z-\tfrac{n-1}{2}\right)\left(zI_{n-1}-C_{i}\right)%
\end{vmatrix}
 \\
&  =\tfrac{P_{A}(z)}{\left(z-\tfrac{n-1}{2}\right)^{n-1}}
\begin{vmatrix}
zI_{n}-B & -J_{n,n-1}\\ 
-nJ_{n-1,n} & \left(z-\tfrac{n-1}{2}\right)\left(zI_{n-1}-C_{i}\right)%
\end{vmatrix}
  \mbox{.}%
\end{align*}
We need to evaluate the determinant
\begin{align*}
&   
\begin{vmatrix}
zI_{n}-B & -J_{n,n-1}\\
-nJ_{n-1,n} & \left(z-\tfrac{n-1}{2}\right)\left(zI_{n-1}-C_{i}\right)
\end{vmatrix}
 \\
&  =n^{n-1} 
\begin{vmatrix}
zI_{n}-B & -J_{n,n-1}\\
-J_{n-1,n} & \tfrac{\left(z-\tfrac{n-1}{2}\right)}{n}\left(zI_{n-1}-C_{i}\right)
\end{vmatrix}
  \mbox{.}%
\end{align*}
Repeating the previous process, we get%
\begin{align*}
&   
\begin{vmatrix}
zI_{n}-B & -J_{n,n-1}\\
-J_{n-1,n} & \tfrac{\left(z-\tfrac{n-1}{2}\right)}{n}\left(zI_{n-1}-C_{i}\right)
\end{vmatrix}
  \\
&  =\tfrac{1}{\left(z-\tfrac{n-1}{2}\right)^{n-1}} 
\begin{vmatrix}
zI_{n}-B & -J_{n,n-1}\\
-\left(z-\tfrac{n-1}{2}\right)J_{n-1,n} & \tfrac{1}{n}\left(z-\tfrac{n-1}{2}\right)^{2}\left(zI_{n-1}-C_{i}\right)
\end{vmatrix}
  \\
&  =\tfrac{1}{\left(  z-\tfrac{n-1}{2}\right)  ^{n-1}}
\begin{vmatrix}
zI_{n}-B & -J_{n,n-1}\\
O_{n-1,n} & \tfrac{1}{n}\left(  z-\tfrac{n-1}{2}\right)  ^{2}\left(zI_{n-1}%
-C_{i}\right)-nJ_{n-1}%
\end{vmatrix}
  \\
&  =\tfrac{1}{\left(  z-\tfrac{n-1}{2}\right)  ^{n-1}}P_{B}(z)\det 
\left( \tfrac{1}{n}(z-\tfrac{n-1}{2})^{2}(zI_{n-1}-C_{i})-nJ_{n-1}\right)  \mbox{.}%
\end{align*}
Hence
\[
P_{M_{i+2n}}(z)=P_{A}(z)P_{B}(z)\det \left(  zI_{n-1}-C_{i}-\tfrac{n^{2}%
}{\left(z-\tfrac{n-1}{2}\right)^{2}}J_{n-1}\right)\mbox{.}
\]
Using the first assertion of Lemma \ref{poly-with indetr}, we have%
\[
P_{M_{i+2n}}(z)=P_{A}(z)P_{B}(z)\left(  \left(  \tfrac{n^{2}}{\left(z-\tfrac{n-1}%
{2}\right)^{2}}+1\right)  P_{C_{i}}(z)-\left(  \tfrac{n^{2}}{\left(z-\tfrac{n-1}{2}\right)^{2}%
}\right)  P_{C_{i}}(-z-1)\right)\mbox{.}
\]
Applying Lemma \ref{regular exp tournament} on the tournament $T_{3}$,
we can write $P_{M_{i+2n}}(z)$ as follows%
\begin{equation}
P_{M_{i+2n}}(z)=F(z)P_{C_{i}}(z)+G(z)\label{expression Mi}\mbox{,}%
\end{equation}
where $F$ and $G$ are rational fractions that depend only on $n$, $P_{A}(z)$,
$P_{B}(z)$, and $P_{C}(z)$.

As $T_{3}$ is not $(n-1)$-spectrally monomorphic, there exists
$j\in [  n]  \setminus  \{  i \}  $ such that $P_{C_{i}%
}(z)\neq P_{C_{j}}(z)$. By \eqref{expression Mi}, we get $P_{M_{i+2n}}(z)\neq
P_{M_{j+2n}}(z)$. Hence, $T$ is not $(3n-1)$-spectrally monomorphic.
\end{proof}

\section{Skew-spectral monomorphy}

We can consider another kind of spectral monomorphy as follows. A
tournament $T$ with adjacency matrix $A$ is $k$-\emph{skew-spectrally
monomorphic} if its skew-adjacency matrix $S=A-A^{\top}$ is $k$-spectrally
monomorphic. The $k$-skew-spectral monomorphy is a weakening of $k$-spectral monomorphy.
Indeed, the following equality \cite[(5.1)]{Gregory1993} shows that the
characteristic polynomial of $A$ determines that of $S$.
\begin{equation}
P_{S}(z)=2^{n-1}\left[  P_{A}\left(  \tfrac{z-1}{2}\right)  +(-1)^{n-1}%
P_{A}\left(  -\tfrac{z+1}{2}\right)  \right]  \label{gregory}%
\end{equation}
 
 Obviously, every tournament is $3$-skew-spectrally monomorphic. Moreover, a tournament $T$ is $k$-skew-spectrally monomorphic if and only if the Hermitian matrix $iS$ is $k$-spectrally monomorphic. 
Then, the results of
\cite{attas2021characterization} provide a characterization of $k$-skew-spectrally
monomorphic $n$-tournaments for $k\in  \{  4,\ldots,n-3 \}  $. In order to state this characterization, we need to define the switching  of tournaments. The \emph{switch} of a tournament $T$, with respect to a subset $X$ of $V$, is
the tournament obtained from $T$ by reversing all the arcs between $X$ and
$V\setminus X$. Two tournaments with the same vertex set are
switching equivalent if and only if their skew-adjacency matrices are
$ \{  \pm1 \}  $-diagonally similar \cite{MOORHOUSE95TW}. It follows
that the switching operation preserves $k$-skew-spectral monomorphy. In particular, every  switch of  a transitive tournament is $k$-skew-spectrally monomorphic. Conversely, using   \cite[Corollary 3.5]{attas2021characterization}, we obtain 
\begin{proposition} 
Let $n$ and $k$ be integers such that $n\geq 8$ and $4\leq k \leq n-4$. An $n$-tournament is $k$-skew-spectrally monomorphic if and only if it is switching equivalent to a transitive tournament. 
\end{proposition} 

In addition to the switch of transitive tournaments, there is another class of $k$-skew-spectrally monomorphic tournaments  for $k=n-3$, which arises
from skew-conference matrices. Recall that a \emph{skew-conference matrix} $C$  is an $n\times n $ skew-symmetric matrix  with $1$ and $-1$ off the diagonal, such that $C^{\top}C=(n-1)I_{n}$. Let $T$ be an  $n$-tournament whose skew-adjacency matrix is a skew-conference matrix. Proposition 4.1 of \cite{attas2021characterization} implies that $T$ is $ k$-skew-spectrally monomorphic for $k\in\{ n-3, n-2,n-1\}$. 
The following proposition provides a complete characterization  of $(n-3)$-skew-spectrally monomorphic $n$-tournaments. It is a consequence of \cite[Theorem 4.2]{attas2021characterization}.
\begin{proposition} 
A tournament with $n\geq 7$ vertices is $(n-3)$-skew-spectrally monomorphic if and only if it is switching equivalent to a transitive tournament or its skew-adjacency matrix is a skew-conference matrix.
\end{proposition} 

For $n\in \{3,4,5\}$, all $n$-tournaments are $(n-2)$-skew-spectrally monomorphic.  For $n\geq 6$, we have three classes of $(n-2)$-skew-spectrally monomorphic $n$-tournaments:

\begin{enumerate}
\item The switches of transitive tournaments.
\item The switches of doubly regular tournaments.
\item Tournaments whose skew-adjacency matrix is a skew-conference matrix.
\end{enumerate}

\begin{problem}
Are there other examples of $(n-2)$-skew-spectrally monomorphic $n$-tournaments?
\end{problem}

As we have seen above, $k$-spectrally monomorphic $n$-tournaments are 
$k$-skew-spectrally monomorphic. 
For every $n\geq4$, we give an example of an $(n-1)$-skew-spectrally monomorphic
$n$-tournament that is not $(n-1)$-spectrally monomorphic.

 Consider the transitive
tournament $T_{n}$ with vertex set $V= \{  v_{1},\ldots,v_{n} \}  $
 such that $v_{i}$ dominates $v_{j}$ if $i<j$. Denote by $R_{n}$ the
tournament obtained from $T_{n}$ by reversing the arc $(  v_{1}%
,v_{n})  $.

\begin{proposition}
The tournament $R_{n}$ is $(n-1)$-skew-spectrally monomorphic
 but not $(n-1)$-spectrally monomorphic.
\end{proposition}

\begin{proof}
 Let $A$ be the adjacency matrix of $R_{n}$ and let
$P_{A}(z)=z^{n}+a_{1}z^{n-1}+\cdots+a_{n}$ be its characteristic polynomial. By
\eqref{eq1}, we have
 $$a_{k}=(-1)^{k}\sum \limits_{\alpha \in \tbinom{\lbrack
n]}{k}}{\det(A[\alpha])}\mbox{.}$$

 Let $k\in \{3,\ldots,n\}$ and let$\  \alpha \in
\tbinom{\lbrack n]}{k}$. Consider the set $Z= \{  v_{i}:i\in
\alpha  \}  $. The adjacency matrix  of the subtournament $R_{n}[Z]$ is $A[\alpha]$.
 If $v_{1}\notin Z$ or $v_{n}\notin Z$ then $R_{n}[Z]$ is
transitive and hence $\det(A[\alpha])=0$. If $ \{  v_{1},v_{n} \}
\subseteq Z$ then $R_{n}[Z]$ is isomorphic to $R_{k}$. By expanding along the
first column, the determinant of the adjacency matrix of $R_{k}$ is
$(-1)^{k+1}$. Then
\[
\sum _{\alpha \in \tbinom{\lbrack n]}{k}}\det(A[\alpha])=(-1)^{k+1}%
\tbinom{n-2}{k-2}\mbox{.}%
\]
Consequently,
\[
P_{A}(z)=z^{n}-{\sum^{n}_{k=3}}\tbinom{n-2}{k-2}z^{n-k}\mbox{.}%
\]

By deleting the vertex $v_{1}$ or $v_{n}$ from the tournament $R_{n}$, we
obtain a transitive tournament with characteristic polynomial $z^{n-1}$. If we
delete any vertex from the set $ \{  v_{2},\ldots,v_{n-1} \}  $, we
obtain a tournament isomorphic to $R_{n-1}$ with characteristic polynomial
$z^{n-1}-{\sum^{n-1}_{k=3}}\tbinom{n-3}{k-2}z^{n-1-k}$. This
proves that $R_{n}$ is not $(n-1)$-spectrally monomorphic.

The switch of the tournament $R_{n}$ with respect to $ \{  v_{1} \}  $ is
transitive. Hence $R_{n}$ is $(n-1)$-skew-spectrally monomorphic.
\end{proof}

For regular tournaments, we have

\begin{proposition}
\label{main3}A regular $n$-tournament is $(n-1)$-spectrally monomorphic
if and only if it is $(n-1)$-skew-spectrally monomorphic.
\end{proposition}

We will prove the following stronger result.

\begin{proposition}
\label{cor A versus S}Let $T$ be a regular $n$-tournament. For $i,j\in  \{
1,\ldots,n \}  $, we have
\[
   P_{A_{j}}(z)-P_{A_{i}}(z)
=\tfrac{  2^{n}z+2^{n-1}}{  2z-n+1}  \left(  P_{S_{j}}(2z+1)-P_{S_{i}}(2z+1)\right)\mbox{,}
\]
where $A$ and $S$ are respectively the adjacency and the skew-adjacency
matrices of $T$.
\end{proposition}

\begin{proof}
From \eqref{gregory}, we get

\begin{align*}
P_{S_{i}}(2z+1)&=2^{n-2}[  P_{A_{i}}(z)+P_{A_{i}}(-z-1)] \mbox{;}\\
P_{S_{j}}(2z+1)&=2^{n-2}[  P_{A_{j}}(z)+P_{A_{j}}(-z-1)]\mbox{.}
\end{align*}
So,
\begin{equation*}
P_{S_{j}}(2z+1)-P_{S_{i}}(2z+1)=2^{n-2}[  P_{A_{j}}(z)-P_{A_{i}%
}(z)+P_{A_{j}}(-z-1)-P_{A_{i}}(-z-1)]\mbox{.}
\end{equation*}
By Lemma \ref{regular exp tournament}, we have
\[
\left(n+2z+1\right)\left(  P_{A_{j}}(z)-P_{A_{i}}(z)\right)  +\left(n-2z-1\right)\left(  P_{A_{j}%
}(-z-1)-P_{A_{i}}(-z-1)\right)  =0
\]
Combining the last two equalities, we get the desired result.
\end{proof}

\bibliography{biblio}

\end{document}